\documentclass[10pt]{article}
\usepackage[utf8]{inputenc}
\usepackage{amsmath}
\usepackage{amssymb}
\usepackage{mathrsfs}
\usepackage{amsfonts}
\usepackage{mathrsfs,amscd,amssymb,amsthm,amsmath,bm,graphicx}
\usepackage{graphicx}
\usepackage{cite}
\usepackage{color}
\usepackage{enumitem}
\makeatletter

\renewcommand{\@seccntformat}[1]{{\csname the#1\endcsname}{\normalsize.}\hspace{.5em}}
\makeatother

\def \[{\begin{equation}}
\def \]{\end{equation}}

\newtheorem{thm}{Theorem}[section]

\begin{document}
\setlength{\baselineskip}{13pt}
\begin{center}{\Large \bf The expected values and limiting behaviours for the Gutman index, Schultz index, multiplicative degree-Kirchhoff index and additive degree-kirchhoff index of a random cyclooctane chain
}

\vspace{4mm}

{\large Jia-Bao Liu $^{1}$, Jiao-Jiao Gu $^{1,*}$, Kang Wang $^{1}$}
\footnotetext{E-mail address: liujiabaoad@163.com, gujiaojiaoajd@163.com, wangkang199804@163.com.\\}
 \footnotetext{* Corresponding author.}
 \vspace{2mm}

{\small $^1$ School of Mathematics and Physics, Anhui Jianzhu University, Hefei 230601, P.R. China\\}
\vspace{2mm}
\end{center}

 {\noindent{\bf Abstract.}\ \
 In this paper, we first introduce the explicit analytical formulas for the expected values of the Gutman and Schultz indices for a random cyclooctane chain $COC_{n}$. Meanwhile, the explicit formulas of the variances of the Gutman and Schultz indices for a random cyclooctane chain are determined and we prove these two indices are asymptotically subject to normal distribution. Furthermore, we are surprised to find the variances of $Kf^{*}(COC_{n})$ and $Kf^{+}(COC_{n})$ for a random cyclooctane chain based on the known results of others' paper and they are asymptotically subject to normal distribution.

 \vspace{2mm}

\section{Introduction} \ \ \ \ \
In this paper, we only take into account simple, undirected and connected graphs. First of all, we review some definitions in graph theory. Suppose $G$ represent a simple undirected graph with $|V_G|=n$ and $|E_G|=m$. The degree of vertex $u$ of the graph $G$ is denoted by $deg(u)$ (or $d(u)$ for short). For more notation, one can be referred to \cite{a1}.

It is well-known that chemical graphs can be used to describe chemical compounds, where vertices correspond to atoms and edges represent the covalent bonds between atoms. An in-depth understanding of the physicochemical properties of chemical compounds is an important part of theoretical chemistry. There are great improvement in predictive methods that link molecular structures to physicochemical properties, and topological index is one of the most important methods.

The distance between vertices $u$ and $v$ of $G$ is the length of a shortest $u$, $v$-path in $G$, denoted by $dis_{G}(u,v)$ (or $d(u,v)$ for short). The famous Wiener index $W(G)$ is the total of distances between all pairs of vertices in the graph $G$, denoted by
\begin{eqnarray*}
W(G)=\sum_{\{u,v\}\subseteq V_{G}}dis_{G}(u,v).
\end{eqnarray*}
This invariant, relating to the distance of graph, was introduced into chemistry in 1947\cite{a2} and mathematics 30 years later\cite{a3}. At present, the Wiener index is a widely used topological index; referred to \cite{a4}.

Now, we introduce to study the weighted Wiener index of a graph. A weighted graph is a graph $G=(V_{G}, E_{G})$ related to the weight function $w$: $V_{G}\longrightarrow {\mathbb{N}}^{+}$, expressed by $(G,w)$. Put $\oplus$ represent one of the arithmetic operations ${+}$,~${-}$,~$\times$,~$\div$, then the weighted Wiener index $W(G,w)$ can be expressed as
\begin{eqnarray}
W(G,w)=\frac{1}{2}\sum_{u \in V_{G}}\sum_{v \in V_{G}}\big(w(u) \oplus w(v)\big)dis_{G}(u,v).
\end{eqnarray}
If $w \equiv 1$ and $\oplus$ stands for the arithmetic operation $\times $, it is straightforward to check that $W(G,w)=W(G)$.

If $\oplus$ stands for the arithmetic operation $\times$ and let $w(\cdot) \equiv deg(\cdot)$, then Eq. (1.1)  is equal to
\begin{eqnarray}
Gut(G)=\frac{1}{2}\sum_{u \in V_{G}}\sum_{v \in V_{G}}\big(deg(u)deg(v)\big)dis_{G}(u,v)=\sum_{\{u,v\}\subseteq V_{G}}\big(deg(u)deg(v)\big)dis_{G}(u,v),
\end{eqnarray}
which is also referred as the Gutman index. For acyclic molecules, $Gut(G)$ has a close relation to $W(G)$ and shows exactly the same structural characteristics of molecules as $W(G)$; see for instance \cite{a5}. Therefore, we should focus more on polycyclic molecules and study the chemical applications and theoretical investigations of the Gutman index.

If $\oplus$ represents the arithmetic operation $+$ and let $w(\cdot) \equiv deg(\cdot)$, Eq. (1.1) is called the Schultz index and equivalent to
\begin{eqnarray}
S(G)=\frac{1}{2}\sum_{u \in V_{G}}\sum_{v \in V_{G}}\big(deg(u)+deg(v)\big)dis_{G}(u,v)=\sum_{\{u,v\}\subseteq V_{G}}\big(deg(u)+deg(v)\big)dis_{G}(u,v).
\end{eqnarray}
 A number of papers are devoted to this graph invariant, referred to \cite{a6,a7}. Information about its chemical applications and properties shall be discovered in \cite{a8}.

For any $u, ~v\in V_{G}$, the effective resistance between $u$ and $v$ can be defined in terms of the potential
difference between $u$ and $v$ when a unit current is maintained from $u$ and $v$, which was on the graph
presented by Klein and Randi\'{c}\cite{a9}, denoted by $r(u, v)$. We can refer to \cite{a10} for details. For non-trees, a famous extension of $W(G)$ is the Kirchhoff index \cite{a9}, denoted by
\begin{eqnarray*}
Kf(G)=\sum_{\{u,v\} \subseteq V_{G}}r(u, v).
\end{eqnarray*}

The multiplicative degree-Kirchhoff index was put forward in 2007 by Chen and Zhang\cite{a11}, which was expressed as
\begin{eqnarray*}
Kf^{*}(G)=\sum_{\{u,v\} \subseteq V_{G}}deg(u)deg(v)r(u, v).
\end{eqnarray*}
This topological index can also be written as
\begin{eqnarray}
Kf^{*}(G)=\sum_{\{u,v\} \subseteq V_{G}}deg(u)deg(v)r(u, v)=2\mid E_{G}\mid \sum_{i=2}^{n}\frac{1}{\lambda_{i}},
\end{eqnarray}
where $0=\lambda_{1}<\lambda_{2}\leq\cdots\leq\lambda_{n}$ are the eigenvalues of $\mathcal{L}(G)$ that stands for the normalized Laplacian matrix, which was introduced by Chung \cite{a1}. In mathematics, chemistry and statistics, their studies have attracted the attention of more and more researchers. A large amount of papers are devoted to them. The readers can refer to \cite{a12} for details.

In 2012, the additive degree-Kirchhoff index was proposed by Gutman, Feng and Yu \cite{a13}, denoted by
\begin{eqnarray}
Kf^{+}(G)=\sum_{\{u,v\} \subseteq V_{G}}\big(deg(u)+deg(v)\big)r(u, v).
\end{eqnarray}
The research on the additive degree-Kirchhoff index is still in its infancy. In practice, the three invariants $Kf(G)$, $Kf^{*}(G)$ and $Kf^{+}(G)$ have a close relationship, and we recommend the recent papers \cite{a14,a15}.

In organic chemistry, chain compounds in particular are an important kind of cycloalkanes. In recent years, a set of macrocyclic aromatic hydrocarbons have attracted the interest of research scholars\cite{a16,a17,a18}, which are known as cyclooctanes and their derivatives. For instance, the synthesis of a number of cyclooctyl pyrazines and quinoxazines has been found by Alamdari et al. \cite{a16} has found. The reader can also refer to \cite{a18} for more information.

The molecular graphs of cyclooctanes are defined as finite 2-connected graphs, whose every interior surface is surrounded by a regular octagon with side length of 1, called octagonal systems\cite{a19,a20}. The tree-like octagonal systems stand for a family of polycyclic conjugated hydrocarbons, and the number of isomers was obtained and showed by the generating functions of tree-like octagonal graphs by Brunvoll et al\cite{a19}. The relationship between the number of perfect matching for a family of octagonal graphs and the Hosoya index of the caterpillar trees was considered by Yang and Zhao\cite{a20}. A lot of interesting combinatorial problems\cite{a21} in octagonal graphs have attracted a great deal of attention from mathematicians. It is defined as a cyclooctane chain if every vertex of the octagonal system is in an octagon and the graph determined by reducing each octagon to a vertex of a octagonal system is a path. There are three kinds of cyclooctane chains at $n=1,~ 2,~3$ as shown in Fig. 1.

\begin{figure}[htbp]
\centering\includegraphics[width=12cm,height=9cm]{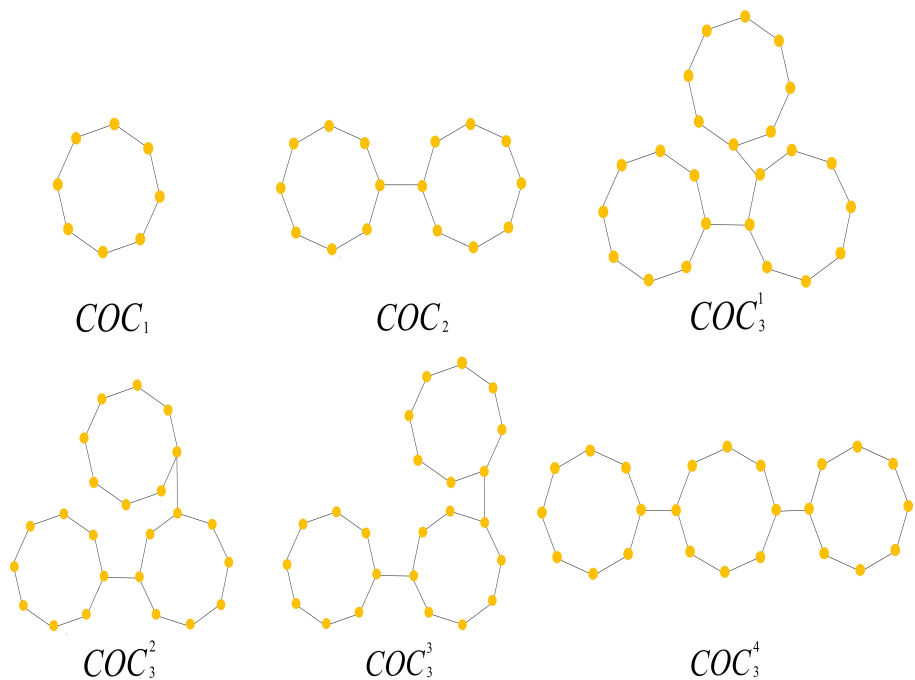}
\caption{ Three kinds of cyclooctane chains.}\label{1}
\end{figure}

\begin{figure}[htbp]
\centering\includegraphics[width=6cm,height=4cm]{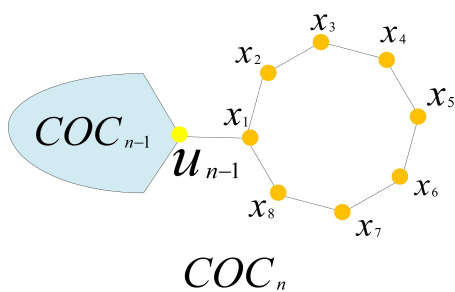}
\caption{  A cyclooctane chain $COC_{n}$.}\label{2}
\end{figure}

In \cite{a22}, a cyclooctane chain $COC_{n}$ with $n$ octagons $O_{1}$, ~$O_{2}$, $\cdot \cdot \cdot$, ~$O_{n}$ is considered as a situation that a new terminal octagon connects by an edge to a cyclooctane chain $COC_{n-1}$ as given in Fig. 2. However, for $n \geq 3$, there are four ways of attaching terminal octagons, and the results can be described as $COC_{n}^{1}$, $COC_{n}^{2}$, $COC_{n}^{3}$ and $COC_{n}^{4}$ as shown in Fig. 3. A random cyclooctane chain, denoted by $COC_{n}(p_{1}, p_{2}, p_{3})$, is a cyclooctane chain given by adding terminal octagons step by step. At each step $k(=3, 4, \cdot \cdot \cdot, n)$, let us choose one of four possible cases at random:

$\bullet$ (i) $COC_{k-1}\longrightarrow COC^{1}_{n}$ with probability $p_{1}$,

$\bullet$ (ii) $COC_{k-1}\longrightarrow COC^{2}_{n}$ with probability $p_{2}$,

$\bullet$ (iii) $COC_{k-1}\longrightarrow COC^{3}_{n}$ with probability $p_{3}$,

$\bullet$ (iiii) $COC_{k-1}\longrightarrow COC^{4}_{n}$ with probability $1-p_{1}-p_{2}-p_{3}$,
\\\text{in which the probabilities $p_{1}$, $p_{2}$ and $p_{3}$ are constants and independent to the step $k$ at the same time.}\par

Motivated by \cite{a23}, we use four random variables $Z^{1}_{n}, Z^{2}_{n},Z^{3}_{n}$ and $Z^{4}_{n}$ to represent our choice. If our choice is $COC^{i}_{n}$, we put $Z^{i}_{n}=1$, otherwise $Z^{i}_{n}=0$ ($i=1,~ 2, ~3, ~4$). One holds that
\begin{eqnarray}
\mathbb{P}(Z^{i}_{n}=1)=P_{i},~~~\mathbb{P}(Z^{i}_{n}=0)=1-P_{i},~~~i=1,~2,~3,~4,
\end{eqnarray}
and $Z^{1}_{n}+Z^{2}_{n}+Z^{3}_{n}+Z^{4}_{n}=1$.

\begin{figure}[htbp]
\centering\includegraphics[width=13cm,height=9cm]{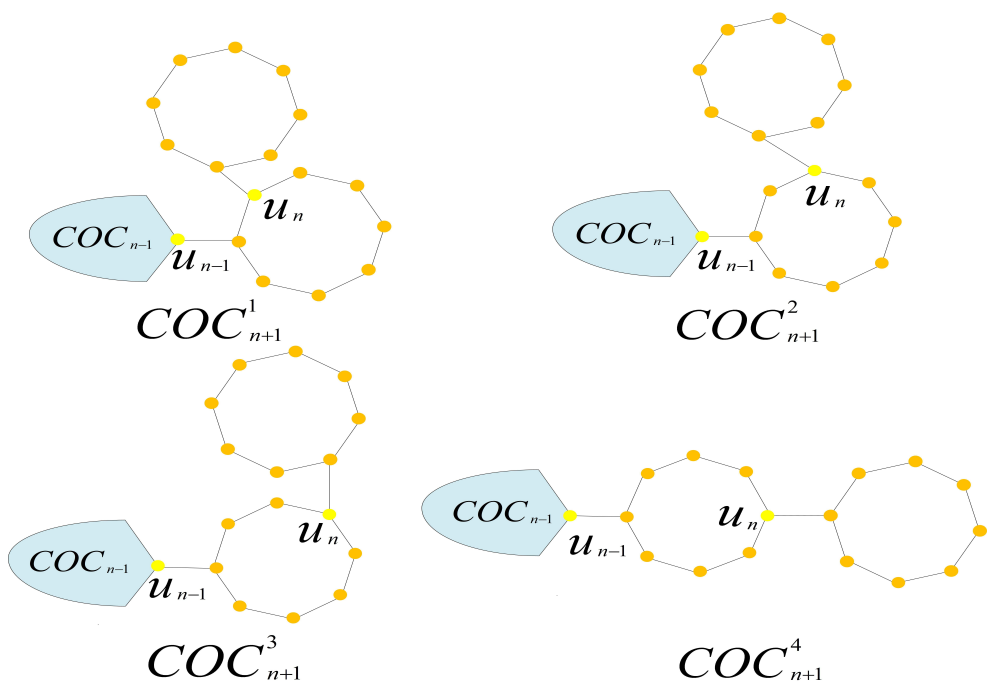}
\caption{ Four attaching ways in cyclooctane chains.}\label{3}
\end{figure}

Yang and Zhang\cite{a24} and Ma et al.\cite{a25} determined explicit formulas of $W(G)$ and $\mathbb{E}\big(W(G)\big)$ for random polyphenylene chains, respectively. Similarly, Huang, Kuang and Deng \cite{a26} obtained $\mathbb{E}\big(Kf(G)\big)$ for random polyphenyl and spiro chains. Wei and Shiu \cite{a27} put forward the expression of $\mathbb{E}\big(W(G)\big)$ for random polygonal chains and proved the asymptotic property of its expected value. At the same way, Zhang, Li, Li and Zhang \cite{a28} obtained the simple formulas of $\mathbb{E}\big(Gut(G_{n})\big)$, $\mathbb{E}\big(S(G_{n})\big)$, $\mathbb{E}\big(Kf^{*}(G_{n})\big)$ and $\mathbb{E}\big(Kf^{+}(G_{n})\big)$ for random polyphenylene chains.

Motivated by \cite{a22,a23,a28,a29}, the rest of the paper is organized as follows. In section 2, we determine the explicit formulas of $\mathbb{E}\big(Gut(COC_{n})\big)$ and $\mathbb{E}\big(S(COC_{n})\big)$ for random cyclooctane chains. In section 3, we obtain explicit formulas of $Var\big(Gut(COC_{n})\big)$ and $Var\big(S(COC_{n})\big)$, and prove that these two indices of random cyclooctane chains asymptotically obey normal distributions. After that, we also establish explicit formulas of $Var\big(Kf^{*}(COC_{n})\big)$ and $Var\big(Kf^{+}(COC_{n})\big)$ based on the known results of others’ paper, and prove that they are asymptotically subject to normal distributions in Section 4.

For the purposes of this paper, we must ensure that the following Hypotheses hold.

\newtheorem{hyp}{Hypothesis}
\begin{hyp}
It is randomly and independently to choice a way attaching the new terminal octagon $O_{n+1}$ to $COC_{n}$, $n=2, 3, \cdot \cdot \cdot$. To be more precisely, the sequences of random variables $\{Z^{1}_{n}, Z^{2}_{n},Z^{3}_{n}, Z^{4}_{n}\}^{\infty}_{n=2}$ are independently and must satisfy Eq. (1.6).
\end{hyp}
\begin{hyp}
For $i\in \{1, 2, 3, 4\}$, we put $0< p_{i} < 1$.
\end{hyp}
Under the condition that Hypotheses 1 and 2,

(a) The analytical expressions of the variances of $Gut(COC_{n})$, $S(COC_{n})$, $Kf^{*}(COC_{n})$ and $Kf^{+}(COC_{n})$ are obtained;

(b) When $n \rightarrow \infty$, we verify that the random variables $Gut(COC_{n})$, $S(COC_{n})$, $Kf^{*}(COC_{n})$ and $Kf^{+}(COC_{n})$ asymptotically obey normal distributions. It is evident to see that
\begin{eqnarray*}
\lim_{n\rightarrow \propto}\sup_{a\in \mathbb{R}}\mid \mathbb{P}\big(\frac{X_{n}-\mathbb{E}(X_{n})}{\sqrt{Var(X_{n})}}\leq a\big)-{\int}^{a}_{-\infty}\frac{1}{\sqrt{2\pi}}e^{-\frac{t^{2}}{2}}dt \mid =0,
\end{eqnarray*}
where $\mathbb{E}(X_{n})$ and $Var(X_{n})$ respectively stand for the expected value and variance of the random variable $X_{n}$.

In this paper, assume that $f(x)$ and $h(x)$ are two functions of $x$. We put $f(x) \asymp h(x)$ if $\lim\limits_{x \rightarrow \infty}\frac{f(x)}{h(x)}=1$, and put $f(x)=O(h(x))$ if $\lim\limits_{x \rightarrow \infty}\frac{f(x)}{h(x)}\leq c$ for $c > 0$.

Please see the following Sections for more details of this paper.

\section{The expected values of $Gut(COC_{n})$ and $S(COC_{n})$ of a random cyclooctane chain}
\ \ \ \
For a random cyclooctane chain $COC_{n}$, we find that $Gut(COC_{n})$ and $S(COC_{n})$ are random variables. Then, we determine the expressions of $\mathbb{E}\big(S(COC_{n})\big)$ and $\mathbb{E}\big(Gut(COC_{n})\big)$ in the section.

In fact, $COC_{n}$ is organized by adding a new terminal octagon $O_{n}$ to $COC_{n-1}$ by an edge, where the vertices of $O_{n}$ are labelled as $x_{1}, ~x_{2},~ x_{3}, ~x_{4},~ x_{5},~ x_{6}, ~x_{7},~ x_{8}$ in clockwise direction. For all $v \in V_{COC_{n-1}}$, one has
\begin{equation}
\begin{split}
&dis(x_{1},v)=dis(u_{n-1},v)+1,~~~dis(x_{2},v)=dis(u_{n-1},v)+2,~~~dis(x_{3},v)=dis(u_{n-1},v)+3,\\
&dis(x_{4},v)=dis(u_{n-1},v)+4,~~~dis(x_{5},v)=dis(u_{n-1},v)+5,~~~dis(x_{6},v)=dis(u_{n-1},v)+4,\\
&dis(x_{7},v)=dis(u_{n-1},v)+3,~~~dis(x_{8},v)=dis(u_{n-1},v)+2,\\
&\sum_{v \in V_{COC_{n-1}}}deg_{COC_{n}}(v)=18n-19, ~~and ~~\sum_{v \in V_{COC_{n}}}deg_{COC_{n+1}}(v)=18n-1.
\end{split}
\end{equation}
In the meantime,
\begin{equation}
\begin{split}
&\sum_{i=1}^{8}deg(x_{i})dis(x_{1},x_{i})=32,~~~\sum_{i=1}^{8}deg(x_{i})dis(x_{2},x_{i})=33,~~~\sum_{i=1}^{8}deg(x_{i})dis(x_{3},x_{i})=34,\\
&\sum_{i=1}^{8}deg(x_{i})dis(x_{4},x_{i})=35,~~~\sum_{i=1}^{8}deg(x_{i})dis(x_{5},x_{i})=36,~~~\sum_{i=1}^{8}deg(x_{i})dis(x_{6},x_{i})=35,\\
&\sum_{i=1}^{8}deg(x_{i})dis(x_{7},x_{i})=34,~~~\sum_{i=1}^{8}deg(x_{i})dis(x_{8},x_{i})=33.\\
\end{split}
\end{equation}

\begin{thm}
For $n\geq 1$, the expression of $\mathbb{E}\big(Gut(COC_{n})\big)$ for a random cyclooctane chain is
\begin{eqnarray*}
\mathbb{E}\big(Gut(COC_{n})\big)&=&(270-162p_{1}-108p_{2}-54p_{3})n^{3}+(486p_{1}+324p_{2}+162p_{3}-90)n^{2}\\
&&+(95-324p_{1}-216p_{2}-108p_{3})n-1.
\end{eqnarray*}
\end{thm}

\noindent\textbf{Proof of Theorem 2.1.}  By Eq. (1.2), one can be convinced that
\begin{eqnarray*}
Gut(COC_{n+1})&=&\sum_{\{u,v\}\subseteq V_{COC_{n}}}deg(u)deg(v)dis(u,v)+\sum_{v\in V_{COC_{n}}}\sum_{x_{i}\in V_{O_{n+1}}}deg(v)deg(x_{i})dis(v,x_{i})\\
&&+\sum_{\{x_{i},x_{j}\}\subseteq V_{O_{n+1}}}deg(x_{i})deg(x_{j})dis(x_{i},x_{j}).
\end{eqnarray*}
Note that
\begin{align*}
&\sum_{\{u,v\}\subseteq V_{COC_{n}}}deg(u)deg(v)dis(u,v)\\
&=\sum_{\{u,v\}\subseteq V_{COC_{n}}\setminus \{u_{n}\}}deg(u)deg(v)dis(u,v)+\sum_{v\in V_{COC_{n}}\setminus \{u_{n}\}}deg_{COC_{n+1}}(u_{n})deg(v)dis(u_{n},v)\\
&=\sum_{\{u,v\}\subseteq V_{COC_{n}}\setminus \{u_{n}\}}deg(u)deg(v)dis(u,v)+\sum_{v\in V_{COC_{n}}\setminus \{u_{n}\}}\big(deg_{COC_{n}}(u_{n})+1\big)deg(v)dis(u_{n},v)\\
&=Gut(COC_{n})+\sum_{v\in V_{COC_{n}}}deg(v)dis(u_{n},v).
\end{align*}
Recall that $deg(x_{1})=3$ and $deg(x_{i})=2$ for $i \in \{2, ~3,~ 4,~ 5, ~6, ~7,~ 8\}$. From Eq. (2.7), we have
\begin{eqnarray*}
\sum_{v\in V_{COC_{n}}}\sum_{x_{i}\in V_{O_{n+1}}}deg(v)deg(x_{i})dis(v,x_{i})&=&\sum_{v\in V_{COC_{n}}}deg(v)\big(\sum_{x_{i}\in V_{O_{n+1}}}deg(x_{i})dis(v,x_{i})\big)\\
&=&\sum_{v\in V_{COC_{n}}}deg(v)\big(17dis(u_{n},v)+49\big)\\
&=&17\sum_{v\in V_{COC_{n}}}deg(v)dis(u_{n},v)+882n-49.
\end{eqnarray*}
From Eq. (2.8), one follows that
\begin{eqnarray*}
\sum_{\{x_{i},x_{j}\}\subseteq V_{O_{n+1}}}deg(x_{i})deg(x_{j})dis(x_{i},x_{j})&=&\frac{1}{2}\sum^{8}_{i=1}deg(x_{i})\big(\sum^{8}_{j=1}deg(x_{j})dis(x_{i},x_{j})\big)\\
&=&\frac{1}{2}(3\times32+2\times 2\times 33+2\times 2\times 34+2\times 2\times 35+2\times 36)\\
&=&288.
\end{eqnarray*}
Then
\begin{eqnarray}
Gut(COC_{n+1})=Gut(COC_{n})+18\sum_{v\in V_{COC_{n}}}deg(v)dis(u_{n},v)+882n+239.
\end{eqnarray}
For a random cyclooctane chain $COC_{n}$, we obtain that $\sum_{v\in V_{COC_{n}}}deg(v)dis(u_{n},v)$ is a random variable. Let
\begin{eqnarray*}
A_{n}:=\mathbb{E}\big(\sum_{v\in V_{COC_{n}}}deg(v)dis(u_{n},v)\big).
\end{eqnarray*}
By using above formula and Eq. (2.9), we can get the following relation for $\mathbb{E}\big(Gut(COC_{n+1})\big)$ of the random cyclooctane chain. One sees
\begin{eqnarray}
\mathbb{E}\big(Gut(COC_{n+1})\big)=\mathbb{E}\big(Gut(COC_{n})\big)+18A_{n}+882n+239.
\end{eqnarray}
Then, we go on to consider the following four possible Cases.

{\bf{Case 1.}} $COC_{n}\longrightarrow COC^{1}_{n+1}$.

In this case, $u_{n}$ (of $COC_{n}$) overlaps with $x_{2}$ or $x_{8}$ (of $O_{n}$). Therefore, $\sum\limits_{v\in V_{COC_{n}}}deg(v)dis(u_{n},v)$ is rewritten as $\sum\limits_{v\in V_{COC_{n}}}deg(v)dis(x_{2},v)$ or $\sum\limits_{v\in V_{COC_{n}}}deg(v)dis(x_{8},v)$ with probability $p_{1}$.

{\bf{Case 2.}} $COC_{n}\longrightarrow COC^{2}_{n+1}$.

In this case, $u_{n}$ (of $COC_{n}$) overlaps with $x_{3}$ or $x_{7}$ (of $O_{n}$). Therefore, $\sum\limits_{v\in V_{COC_{n}}}deg(v)dis(u_{n},v)$ is rewritten as $\sum\limits_{v\in V_{COC_{n}}}deg(v)dis(x_{3},v)$ or $\sum\limits_{v\in V_{COC_{n}}}deg(v)dis(x_{7},v)$ with probability $p_{2}$.

{\bf{Case 3.}} $COC_{n}\longrightarrow COC^{3}_{n+1}$.

 In this case, $u_{n}$ (of $COC_{n}$) overlaps with $x_{4}$ or $x_{6}$ (of $O_{n}$). Therefore, $\sum\limits_{v\in V_{COC_{n}}}deg(v)dis(u_{n},v)$ is rewritten as $\sum\limits_{v\in V_{COC_{n}}}deg(v)dis(x_{4},v)$ or $\sum\limits_{v\in V_{COC_{n}}}deg(v)dis(x_{6},v)$ with probability $p_{3}$.

{\bf{Case 4.}} $COC_{n}\longrightarrow COC^{4}_{n+1}$.

In this case, $u_{n}$ (of $COC_{n}$) overlaps with $x_{5}$ (of $O_{n}$). Therefore, $\sum\limits_{v\in V_{COC_{n}}}deg(v)dis(u_{n},v)$ is rewritten as $\sum\limits_{v\in V_{COC_{n}}}deg(v)dis(x_{5},v)$ with probability $1-p_{1}-p_{2}-p_{3}$.

Together with the above Cases, by applying the expectation operator and Eq. (2.7), one follows that
\begin{eqnarray*}
A_{n}&=&p_{1}\sum\limits_{v\in V_{COC_{n}}}deg(v)dis(x_{2},v)+p_{2}\sum\limits_{v\in V_{COC_{n}}}deg(v)dis(x_{3},v)+p_{3}\sum\limits_{v\in V_{COC_{n}}}deg(v)dis(x_{4},v)\\
&&+(1-p_{1}-p_{2}-p_{3})\sum\limits_{v\in V_{COC_{n}}}deg(v)dis(x_{5},v),
\end{eqnarray*}
where
\begin{eqnarray*}
&\sum\limits_{v\in V_{COC_{n}}}deg(v)dis(x_{2},v)=\sum\limits_{v\in V_{COC_{n-1}}}deg(v)dis(u_{n-1},v)+2\sum\limits_{v\in V_{COC_{n-1}}}deg(v)+33,\\
&\sum\limits_{v\in V_{COC_{n}}}deg(v)dis(x_{3},v)=\sum\limits_{v\in V_{COC_{n-1}}}deg(v)dis(u_{n-1},v)+3\sum\limits_{v\in V_{COC_{n-1}}}deg(v)+34,\\
&\sum\limits_{v\in V_{COC_{n}}}deg(v)dis(x_{4},v)=\sum\limits_{v\in V_{COC_{n-1}}}deg(v)dis(u_{n-1},v)+4\sum\limits_{v\in V_{COC_{n-1}}}deg(v)+35,\\
and &\sum\limits_{v\in V_{COC_{n}}}deg(v)dis(x_{5},v)=\sum\limits_{v\in V_{COC_{n-1}}}deg(v)dis(u_{n-1},v)+5\sum\limits_{v\in V_{COC_{n-1}}}deg(v)+36.
\end{eqnarray*}
Then, we obtain
\begin{eqnarray*}
A_{n}&=&p_{1}\big[A_{n-1}+2(18n-19)+33\big]+p_{2}\big[A_{n-1}+3(18n-19)+34\big]+p_{3}\big[A_{n-1}+4(18n-19)+35\big]\\
&&+(1-p_{1}-p_{2}-p_{3})\big[A_{n-1}+5(18n-19)+36\big]\\
&=&A_{n-1}+(90-54p_{1}-36p_{2}-18p_{3})n+(54p_{1}+36p_{2}+18p_{3}-59).
\end{eqnarray*}
Meanwhile, for $n=1$, the boundary condition is
\begin{eqnarray*}
A_{1}=\sum_{v\in V_{COC_{1}}}deg(v)dis(u_{1},v)=32.
\end{eqnarray*}
Using above condition and the recurrence relation with respect to $A_{n}$, it is no hard to obtain
\begin{eqnarray}
A_{n}&=&(45-27p_{1}-18p_{2}-9p_{3})n^{2}+(27p_{1}+18p_{2}+9p_{3}-14)n+1.
\end{eqnarray}
From Eq. (2.10), it holds that
\begin{eqnarray*}
\mathbb{E}\big(Gut(COC_{n+1})\big)&=&\mathbb{E}(Gut(COC_{n}))+18A_{n}+882n+239\\
&=&\mathbb{E}(Gut(COC_{n}))+18\bigg((45-27p_{1}-18p_{2}-9p_{3})n^{2}+(27p_{1}+18p_{2}+9p_{3}-14)n\\
&&+1\bigg)+882n+239.
\end{eqnarray*}
For $n=1$, we obtain $\mathbb{E}\big(Gut(COC_{1})\big)=256$.

Similarly, according to the recurrence relation related to $\mathbb{E}\big(Gut(COC_{n})\big)$, we have
\begin{eqnarray*}
\mathbb{E}\big(Gut(COC_{n})\big)&=&(270-162p_{1}-108p_{2}-54p_{3})n^{3}+(486p_{1}+324p_{2}+162p_{3}-90)n^{2}\\
&&+(95-324p_{1}-216p_{2}-108p_{3})n-1,
\end{eqnarray*}
as desired.\hfill\rule{1ex}{1ex}\

\begin{thm}
The expression of $\mathbb{E}\big(S(COC_{n})\big)$ for a random cyclooctane chain is
\begin{eqnarray*}
\mathbb{E}\big(S(COC_{n})\big)&=&(240-144p_{1}-96p_{2}-48p_{3})n^{3}+(432p_{1}+288p_{2}+144p_{3}-80)n^{2}\\
&&+(96-288p_{1}-192p_{2}-96p_{3})n.
\end{eqnarray*}
\end{thm}

\noindent\textbf{Proof of Theorem 2.2.} Notice that a cyclooctane chain $COC_{n+1}$ is organized by adding a new terminal octagon $O_{n+1}$ to $COC_{n}$ by an edge. By Eq. (1.3), one has
\begin{eqnarray}
S(COC_{n+1})=\Lambda_{1}+\Lambda_{2}+\Lambda_{3},
\end{eqnarray}
where
\begin{eqnarray*}
\Lambda_{1}=\sum_{\{u,v\}\subseteq V_{COC_{n}}}(deg(u)+deg(v))dis(u,v),\\
\Lambda_{2}=\sum_{v\in V_{COC_{n}}}\sum_{x_{i}\in V_{O_{n+1}}}(deg(v)+deg(x_{i}))dis(v,x_{i}),\\
\Lambda_{3}=\sum_{\{x_{i},x_{j}\}\subseteq V_{O_{n+1}}}(deg(x_{i})+deg(x_{j}))dis(x_{i},x_{j}).\\
\end{eqnarray*}

Put $dis(u_{n}\mid COC_{n}):=\sum\limits_{v\in V_{COC_{n}}}dis(u_{n},v)$. It is routine to check that
\begin{eqnarray*}
\Lambda_{1}&=&\sum_{\{u,v\}\subseteq V_{COC_{n}}\setminus \{u_{n}\}}(deg(u)+deg(v))dis(u,v)+\sum_{v\in V_{COC_{n}}\setminus \{u_{n}\}}(deg_{COC_{n+1}}(u_{n})+deg(v))dis(u_{n},v)\\
&=&\sum_{\{u,v\}\subseteq V_{COC_{n}}\setminus \{u_{n}\}}(deg(u)+deg(v))dis(u,v)+\sum_{v\in V_{COC_{n}}\setminus \{u_{n}\}}\big(deg_{COC_{n}}(u_{n})+1+deg(v)\big)dis(u_{n},v)\\
&=&S(COC_{n})+dis(u_{n}\mid COC_{n}).
\end{eqnarray*}
Note that $COC_{n}$ has $8n$ vertices. Note that $deg(x_{1})=3$ and $deg(x_{i})=2$ for $i \in \{2,~ 3, ~4, ~5, ~6, ~7, ~8\}$. By Eq. (2.7), we have
\begin{eqnarray*}
\Lambda_{2}&=&\sum_{v\in V_{COC_{n}}}\sum_{x_{i}\in V_{O_{n+1}}}deg(v)dis(v,x_{i})+\sum_{v\in V_{COC_{n}}}\sum_{x_{i}\in V_{O_{n+1}}}deg(x_{i})dis(x_{i},v)\\
&=&\sum_{v\in V_{COC_{n}}}deg(v)\big(\sum_{x_{i}\in V_{O_{n+1}}}dis(v,x_{i})\big)+\sum_{v\in V_{COC_{n}}}\big(\sum_{x_{i}\in V_{O_{n+1}}}deg(x_{i})dis(x_{i},v)\big)\\
&=&8\sum_{v\in V_{COC_{n}}}deg(v)dis(u_{n},v)+24(18n-1)+17\sum_{v\in V_{COC_{n}}}dis(u_{n},v)+49\times 8n\\
&=&8\sum_{v\in V_{COC_{n}}}deg(v)dis(u_{n},v)+17dis(u_{n}\mid COC_{n})+824n-24.
\end{eqnarray*}
Note that $\sum\limits_{i=1}^{8}dis(x_{k},x_{i})=16$ for $k=1, ~2,\cdot \cdot \cdot, ~8$. From Eq. (2.8), one sees that
\begin{eqnarray*}
\Lambda_{3}&=&\sum_{\{x_{i},x_{j}\}\subseteq V_{O_{n+1}}}(deg(x_{i})+deg(x_{j}))dis(x_{i},x_{j})=\frac{1}{2}\sum^{8}_{i=1}\sum^{8}_{j=1}(deg(x_{i})+deg(x_{j}))dis(x_{i},x_{j})\\
&=&\sum^{8}_{i=1}\sum^{8}_{j=1}deg(x_{i})dis(x_{i},x_{j})=16(3+2\times 7)=272.
\end{eqnarray*}
Then, Eq. (2.12) can be rewritten as
\begin{eqnarray}
S(COC_{n+1})=S(COC_{n})+8\sum_{v\in V_{COC_{n}}}deg(v)dis(u_{n},v)+18dis(u_{n}\mid COC_{n})+824n+248.
\end{eqnarray}
For a random cyclooctane chain $COC_{n}$, we know that $dis(u_{n}\mid COC_{n})$ is a random variable. Let
\begin{eqnarray*}
B_{n}:=\mathbb{E}\big(dis(u_{n}\mid COC_{n})\big).
\end{eqnarray*}
According to above formula and Eq. (2.13), we have the following relation for $\mathbb{E}\big(S(COC_{n+1})\big)$ of a random cyclooctane chain. It holds that
\begin{eqnarray}
\mathbb{E}\big(S(COC_{n+1})\big)=\mathbb{E}\big(S(COC_{n})\big)+8A_{n}+18B_{n}+824n+248.
\end{eqnarray}
We proceed by taking into account the following four Cases.

{\bf{Case 1.}} $COC_{n}\longrightarrow COC^{1}_{n+1}$.

In this case, $u_{n}$ (of $COC_{n}$) coincides with $x_{2}$ or $x_{8}$ (of $O_{n}$). Then, $dis(u_{n}\mid COC_{n})$ is given by $dis(x_{2}\mid COC_{n})$ or $dis(x_{8}\mid COC_{n})$ with probability $p_{1}$.

{\bf{Case 2.}} $COC_{n}\longrightarrow COC^{2}_{n+1}$.

In this case, $u_{n}$ (of $COC_{n}$) coincides with $x_{3}$ or $x_{7}$ (of $O_{n}$). Then, $dis(u_{n}\mid COC_{n})$ is given by $dis(x_{3}\mid COC_{n})$ or $dis(x_{7}\mid COC_{n})$ with probability $p_{2}$.

{\bf{Case 3.}} $COC_{n}\longrightarrow COC^{3}_{n+1}$.

In this case, $u_{n}$ (of $COC_{n}$) coincides with $x_{4}$ or $x_{6}$ (of $O_{n}$). Then, $dis(u_{n}\mid COC_{n})$ is given by $dis(x_{4}\mid COC_{n})$ or $dis(x_{6}\mid COC_{n})$ with probability $p_{3}$.

{\bf{Case 4.}} $COC_{n}\longrightarrow COC^{4}_{n+1}$.

 In this case, $u_{n}$ (of $COC_{n}$) coincides with $x_{5}$ (of $O_{n}$). Then, $dis(u_{n}\mid COC_{n})$ is given by $dis(x_{5}\mid COC_{n})$ with probability $1-p_{1}-p_{2}-p_{3}$.

Together with Cases 1, 2, 3 and 4, we obtain $B_{n}$ equals to
\begin{align*}
p_{1}\cdot dis(x_{2}\mid COC_{n})+p_{2}\cdot dis(x_{3}\mid COC_{n})+p_{3}\cdot dis(x_{4}\mid COC_{n})+(1-p_{1}-p_{2}-p_{3})\cdot dis(x_{5}\mid COC_{n}),
\end{align*}
where
\begin{eqnarray*}
&dis(x_{2}\mid COC_{n})=dis(u_{n-1}\mid COC_{n-1})+16n,\\
&dis(x_{3}\mid COC_{n})=dis(u_{n-1}\mid COC_{n-1})+24n-8,\\
&dis(x_{4}\mid COC_{n})=dis(u_{n-1}\mid COC_{n-1})+32n-16,\\
&dis(x_{5}\mid COC_{n})=dis(u_{n-1}\mid COC_{n-1})+40n-24.
\end{eqnarray*}
Then, as an immediate consequence, we have
\begin{eqnarray*}
B_{n}&=&p_{1}\big(B_{n-1}+16n\big)+p_{2}\big(B_{n-1}+24n-8\big)+p_{3}\big(B_{n-1}+32n-16\big)\\
&&+(1-p_{1}-p_{2}-p_{3})\big(B_{n-1}+40n-24\big)\\
&=&B_{n-1}+(40-24p_{1}-16p_{2}-8p_{3})n+(24p_{1}+16p_{2}+8p_{3}-24).
\end{eqnarray*}
Put $n=1$, we can easily find that
\begin{eqnarray*}
B_{1}=\mathbb{E}(dis(u_{1}\mid COC_{1}))=16.
\end{eqnarray*}
Then, using above formula and the recurrence relation, it is routine to check that
\begin{eqnarray*}
B_{n}&=&(20-12p_{1}-8p_{2}-4p_{3})n^{2}+(12p_{1}+8p_{2}+4p_{3}-4)n.
\end{eqnarray*}
From Eq. (2.11), we have
\begin{eqnarray*}
A_{n}&=&(45-27p_{1}-18p_{2}-9p_{3})n^{2}+(27p_{1}+18p_{2}+9p_{3}-14)n+1.
\end{eqnarray*}
By using Eq. (2.14), we get
\begin{eqnarray*}
\mathbb{E}\big(S(COC_{n+1})\big)&=&\mathbb{E}(S(COC_{n}))+8A_{n}+18B_{n}+824n+248\\
&=&\mathbb{E}(S(COC_{n}))+8\bigg((45-27p_{1}-18p_{2}-9p_{3})n^{2}+(27p_{1}+18p_{2}+9p_{3}-14)n+1\bigg)\\
&&+18\bigg((20-12p_{1}-8p_{2}-4p_{3})n^{2}+(12p_{1}+8p_{2}+4p_{3}-4)n\bigg)+824n+248.
\end{eqnarray*}
As an immediate consequence, we find $\mathbb{E}(S(COC_{1}))=256$. Then, we arrive at
\begin{eqnarray*}
\mathbb{E}(S(COC_{n}))&=&(240-144p_{1}-96p_{2}-48p_{3})n^{3}+(432p_{1}+288p_{2}+144p_{3}-40)n^{2}\\
&&+(56-288p_{1}-192p_{2}-96p_{3})n,
\end{eqnarray*}
as desired.\hfill\rule{1ex}{1ex}\

\section{The limiting behaviours for $Gut(COC_{n})$ and $S(COC_{n})$ of a random cyclooctane chain}
In Section 3, we obtain the explicit analytical expressions for $Var\big(Gut(COC_{n})\big)$ and $Var\big(S(COC_{n})\big)$. For random cyclooctane chains, we prove the Gutman and Schultz indices asymptotically obey normal distributions. We use the same notation as those used at Section 2.

\begin{thm}
Suppose Hypotheses 1 and 2 are true, then the next results hold. $(i)$ The variance of the Gutman index is denoted by
\begin{eqnarray*}
Var\big(Gut(COC_{n})\big)&=&\frac{1}{30}\bigg({\sigma}^{2}_{1}n^{5}-5r_{1}n^{4}+10\tilde{\sigma}^{2}_{1}n^{3}+(65r_{1}-30{\sigma}^{2}_{1}-45\tilde{\sigma}^{2}_{1})n^{2}\\
&&+(59{\sigma}^{2}_{1}+65\tilde{\sigma}^{2}_{1}-120r_{1})n+(60r_{1}-30{\sigma}^{2}_{1}-30\tilde{\sigma}^{2}_{1})\bigg),
\end{eqnarray*}
where
\begin{align*}
\sigma^{2}_{1}&=648^{2}p_{1}+972^{2}p_{2}+1268^{2}p_{3}+1620^{2}(1-p_{1}-p_{2}-p_{3})\\
&~~~~-\big(648p_{1}+972p_{2}+1268p_{3}+1620(1-p_{1}-p_{2}-p_{3})\big)^{2},\\
\tilde{\sigma}^{2}_{1}&=90^{2}p_{1}+414^{2}p_{2}+738^{2}p_{3}+1062^{2}(1-p_{1}-p_{2}-p_{3})\\
&~~~~-\big(90p_{1}+414p_{2}+738p_{3}+1062(1-p_{1}-p_{2}-p_{3})\big)^{2},\\
r_{1}&=-\big(648p_{1}+972p_{2}+1268p_{3}+1620(1-p_{1}-p_{2}-p_{3})\big)\big(90p_{1}+414p_{2}+738p_{3}+1062(1-p_{1}-p_{2}-p_{3})\big)\\
&~~~~+648\cdot 90p_{1}+972\cdot 414p_{2}+1268\cdot 738p_{3}+1620\cdot 1062(1-p_{1}-p_{2}-p_{3}).
\end{align*}
(ii) For $n\rightarrow \infty$, $Gut(COC_{n})$ asymptotically obeys normal distributions. One has
\begin{eqnarray*}
\lim_{n\rightarrow \propto}\sup_{a\in \mathbb{R}}\mid \mathbb{P}\big(\frac{Gut(COC_{n})-\mathbb{E}\big(Gut(COC_{n})\big)}{\sqrt{Var\big(Gut(COC_{n})\big)}}\leq a\big)-{\int}^{a}_{-\infty}\frac{1}{\sqrt{2\pi}}e^{-\frac{t^{2}}{2}}dt \mid =0.
\end{eqnarray*}
\end{thm}

\noindent\textbf{Proof of Theorem 3.1.} Let $C_{n}=18\sum\limits_{v \in V_{COC_{n}}}deg(v)dis(u_{n},v)$. Then by Eq. (2.9), we obtain
\begin{eqnarray}
Gut(COC_{n+1})=Gut(COC_{n})+C_{n}+882n+239.
\end{eqnarray}
Recalling that $Z_{n}^{1}$, $Z_{n}^{2}$, $Z_{n}^{3}$ and $Z_{n}^{4}$ are random variables which stand for our choice to construct $COC_{n+1}$ by $COC_{n}$. We have the next four Facts.

\noindent\textbf{\bf Fact 3.1.1.} $C_{n}Z_{n}^{1}=(C_{n-1}+648n-90)Z_{n}^{1}.$

\noindent\textbf{\bf Proof.} If $Z_{n}^{1}$ = 0, the result is obvious. Then, we only take into account $Z_{n}^{1}$ = 1, which implies $COC_{n} \rightarrow COC^{1}_{n+1}$. In this case, $u_{n}$ (of $COC_{n}$) overlaps with $x_{2}$ or $x_{8}$ (of $O_{n}$), see Fig. 3. In this situation, by using Eq. (2.7) - Eq. (2.8), $C_{n}$ becomes
\begin{eqnarray*}
18\sum_{v \in V_{COC_{n}}}d(v)d(x_{2},v)&=&18\sum_{v \in V_{COC_{n-1}}}deg(v)dis(x_{2},v)+18\sum_{v \in V_{O_{n}}}deg(v)dis(x_{2},v)\\
&=&18\sum_{v \in V_{COC_{n-1}}}deg(v)\big(dis(v,u_{n-1})+dis(x_{2},u_{n-1})\big)+18\times 33\\
&=&18\sum_{v \in V_{COC_{n-1}}}deg(v)\big(dis(v,u_{n-1})+2\big)+18\times 33\\
&=&C_{n-1}+648n-90.
\end{eqnarray*}
Thus, we conclude the desired Fact.

\noindent\textbf{\bf Fact 3.1.2.} $C_{n}Z_{n}^{2}=(C_{n-1}+972n-414)Z_{n}^{2}.$

Similar to the proof of Fact 3.1.1, we only consider the fact $Z_{n}^{2}=1$, that is $COC_{n} \rightarrow COC^{2}_{n+1}$. In the same way, we omit the details.

\noindent\textbf{\bf Fact 3.1.3.} $C_{n}Z_{n}^{3}=(C_{n-1}+1296n-738)Z_{n}^{3}.$

We only consider the fact $Z_{n}^{3}=1$, that is $COC_{n} \rightarrow COC^{3}_{n+1}$. The proof is also similar to the above facts and we omit the details.

\noindent\textbf{\bf Fact 3.1.4.} $C_{n}Z_{n}^{4}=(C_{n-1}+1620n-1062)Z_{n}^{4}.$

Considering the fact $Z_{n}^{4}=1$ which is $COC_{n} \rightarrow COC^{4}_{n+1}$, the details are omitted here.

Noting that $Z_{n}^{1}+Z_{n}^{2}+Z_{n}^{3}+Z_{n}^{4}=1$, by the above discussions, it holds that
\begin{eqnarray*}
C_{n}&=&C_{n}(Z_{n}^{1}+Z_{n}^{2}+Z_{n}^{3}+Z_{n}^{4})\\
&=&(C_{n-1}+648n-90)Z_{n}^{1}+(C_{n-1}+972n-414)Z_{n}^{2}+(C_{n-1}+1296n-738)Z_{n}^{3}\\
&&+(C_{n-1}+1620n-1062)Z_{n}^{4}\\
&=&C_{n-1}+(648Z_{n}^{1}+972Z_{n}^{2}+1296Z_{n}^{3}+1620Z_{n}^{4})n-(90Z_{n}^{1}+414Z_{n}^{2}+738Z_{n}^{3}+1062Z_{n}^{4})\\
&=&C_{n-1}+nU_{n}-V_{n},
\end{eqnarray*}
where for each $n$,

~~~~$U_{n}=648Z_{n}^{1}+972Z_{n}^{2}+1296Z_{n}^{3}+1620Z_{n}^{4}$,~~~~$V_{n}=90Z_{n}^{1}+414Z_{n}^{2}+738Z_{n}^{3}+1062Z_{n}^{4}$.
\\\text{Therefore, by Eq. (3.15), it follows that}\par

\begin{align}
Gut(COC_{n})=& Gut(COC_{1})+\sum_{a=1}^{n-1}C_{a}+\sum_{a=1}^{n-1}(882a+239)\nonumber\\
=& Gut(COC_{1})+\sum_{a=1}^{n-1}(\sum_{b=1}^{a-1}(C_{b+1}-C_{b})+C_{1})+\sum_{a=1}^{n-1}(882a+239)\nonumber\\
=& Gut(COC_{1})+\sum_{a=1}^{n-1}\sum_{b=1}^{a-1}(C_{b+1}-C_{b})+(n-1)C_{1}+\sum_{a=1}^{n-1}(882a+239)\nonumber\\
=& Gut(COC_{1})+\sum_{a=1}^{n-1}\sum_{b=1}^{a-1}((b+1)U_{b+1}-V_{b+1})+O(n^2).
\end{align}
By direct calculation, we put

~~~~~~~~~~$Var(U_{b})=\sigma^{2}_{1}$,~~~~$Var(V_{b})=\tilde{\sigma}^{2}_{1}$,~~~~$Cov(U_{b},V_{b})=r_{1}$,
\\\text{where for any two random variables $X$, $Y$, $Cov(X,Y):=\mathbb{E}(XY)-\mathbb{E}(X)\mathbb{E}(Y)$.}\par

By the properties of variance, Eq. (3.16) and interchanging the orders of $a$ and $b$, it is straightforward to check that
\begin{align*}
&Var\big(Gut(COC_{n})\big)\\
&=Var\big(\sum_{a=1}^{n-1}\sum_{b=1}^{a-1}(b+1)U_{b+1}-V_{b+1}\big)=Var\big(\sum_{b=1}^{n-2}\sum_{a=b+1}^{n-1}(b+1)U_{b+1}-V_{b+1}\big)\\
&=Var\big(\sum_{b=1}^{n-2}((b+1)U_{b+1}-V_{b+1})(n-b-1)\big)=\sum_{b=1}^{n-2}(n-b-1)^{2}Var\big((b+1)U_{b+1}-V_{b+1}\big)\\
&=\sum_{b=1}^{n-2}(n-b-1)^{2}Cov\big((b+1)U_{b+1}-V_{b+1},(b+1)U_{b+1}-V_{b+1}\big)\\
&=\sum_{b=1}^{n-2}(n-b-1)^{2}\bigg((b+1)^2Cov(U_{b+1},U_{b+1})-2(b+1)Cov(U_{b+1},V_{b+1})+Cov(V_{b+1},V_{b+1})\bigg)\\
&=\sum_{b=1}^{n-2}(n-b-1)^{2}\big((b+1)^2\sigma^{2}_{1}-2(b+1)r_{1}+\tilde{\sigma}^{2}_{1}\big).\\
\end{align*}
By using a computer, the above expression indicates the result $Theorem ~3.1.~~ (i)$.

Now we will go on to the proof of $Theorem ~3.1.~~ (ii)$. Firstly, for any $n\in \mathbb{N}$, let

~~~~~$\mathcal{U}_{n}=\sum\limits_{a=1}^{n-1}\sum\limits_{b=1}^{a-1}(b+1)U_{b+1}$,~~~~$\mathcal{V}_{n}=\sum\limits_{a=1}^{n-1}\sum\limits_{b=1}^{a-1}V_{b+1}$, ~~$\mu= \mathbb{E}(U_{b})$, ~~and~~   $\phi (t)=\mathbb{E}(e^{t(U_{b}-\mu)})$.
\\\text{By these notations, obviously, we have}\par

$e^{t(\mathcal{U}_{n}-\mathbb{E}(\mathcal{U}_{n}))}=e^{t\sum\limits_{a=1}^{n-1}\sum\limits_{b=1}^{a-1}(b+1)(U_{b+1}-\mu)}=e^{t\sum\limits_{b=1}^{n-2}\sum\limits_{a=b+1}^{n-1}(b+1)(U_{b+1}-\mu)}=e^{t\sum\limits_{b=1}^{n-2}(n-b-1)(b+1)(U_{b+1}-\mu)}$,
\\\text{then}\par
\begin{align}
\mathbb{E}\big(e^{t(\mathcal{U}_{n}-\mathbb{E}(\mathcal{U}_{n}))}\big)=&\mathbb{E}\big(e^{t\sum\limits_{b=1}^{n-2}(n-b-1)(b+1)(U_{b+1}-\mu)}\big)=\prod\limits_{b=1}^{n-2}\mathbb{E}\big(e^{t(n-b-1)(b+1)(U_{b+1}-\mu)}\big)\nonumber\\
=&\prod\limits_{b=1}^{n-2}\phi \big(t(n-b-1)(b+1)\big),
\end{align}
and for $k>0$,
\begin{eqnarray}
\mathcal{V}_{n}=\sum\limits_{a=1}^{n-1}\sum\limits_{b=1}^{a-1}V_{b+1}\leq kn^2.
\end{eqnarray}
Noting that

~~$Var(Gut(COC_{n})) \asymp \frac{1}{30}\sigma_{1}^2n^5$,~~ $\phi (t)=1+\frac{\sigma_{1}^2}{2}t^2+O(t^2)$, ~~and ~~$\sum\limits^{n-2}_{b=1}(b+1)^2(n-b-1)^2\asymp \frac{n^5}{30}$.
\\\text{By Taylor's formula and Eqs. (3.16)-(3.18), one holds that}\par
\begin{align*}
&\lim\limits_{n\rightarrow \infty}\mathbb{E}\exp\bigg\{t\frac{Gut(COC_{n})-\mathbb{E}(Gut(COC_{n}))}{\sqrt{Var(Gut(COC_{n}))}}\bigg\}\\
&=\lim\limits_{n\rightarrow \infty}\mathbb{E}\exp\bigg\{t\frac{\big(Gut(COC_{1})+\mathcal{U}_{n}-\mathcal{V}_{n}+O(n^2)\big)-\mathbb{E}\big(Gut(COC_{1})+\mathcal{U}_{n}-\mathcal{V}_{n}+O(n^2)\big)}{\frac{\sigma_{1} n^{\frac{5}{2}}}{\sqrt{30}}}\bigg\}\\
&=\lim\limits_{n\rightarrow \infty}\mathbb{E}\exp\bigg\{t\frac{\sqrt{30}\big(\mathcal{U}_{n}-\mathbb{E}(\mathcal{U}_{n})\big)}{\sigma_{1} n^{\frac{5}{2}}}\bigg\}=\lim\limits_{n\rightarrow \infty}\prod\limits_{b=1}^{n}\phi \big(\frac{\sqrt{30}t(b+1)(n-b-1)}{\sigma_{1} n^{\frac{5}{2}}}\big)\\
&=\lim\limits_{n\rightarrow \infty}\exp\bigg\{\sum\limits_{b=1}^{n}\ln\phi \big(\frac{\sqrt{30}t(b+1)(n-b-1)}{\sigma_{1} n^{\frac{5}{2}}}\big)\bigg\}\\
&=\lim\limits_{n\rightarrow \infty}\exp\bigg\{\sum\limits_{b=1}^{n}\ln\big(1+\frac{\sigma_{1}^2}{2}\frac{30t^2(b+1)^2(n-b-1)^2}{\sigma_{1}^2n^5}+O(\frac{1}{n})\big)\bigg\}\\
&=\lim\limits_{n\rightarrow \infty}\exp\bigg\{\sum\limits_{b=1}^{n-2}\big(\frac{\sigma_{1}^2}{2}\frac{30t^2(b+1)^2(n-b-1)^2}{\sigma_{1}^2n^5}+O(\frac{1}{n})\big)\bigg\}\\
&=e^{\frac{t^2}{2}},
\end{align*}
 Assume that $\mathbb{I}$ is a complex number with $\mathbb{I}^2=-1$. We use $\mathbb{I}t$ instead of $t$ and one has
\begin{align*}
\lim\limits_{n\rightarrow \infty}\mathbb{E}\exp\bigg\{\mathbb{I}t\frac{Gut(COC_{n})-\mathbb{E}(Gut(COC_{n}))}{\sqrt{Var(Gut(COC_{n}))}}\bigg\}=e^{-\frac{t^2}{2}}.
\end{align*}
According to the above formula (\cite{a30}, Chapter 1) and the theory of continuity of probability characteristic functions (\cite{a31}, Chapter 15), we complete the proof of $Theorem ~3.1.~~ (ii)$.\hfill\rule{1ex}{1ex}\

Now we consider $S(COC_{n})$. In $Theorem ~~2.2$, we prove that
\begin{eqnarray*}
\mathbb{E}(S(COC_{n}))&=&(240-144p_{1}-96p_{2}-48p_{3})n^{3}+(432p_{1}+288p_{2}+144p_{3}-80)n^{2}\\
&&+(96-288p_{1}-192p_{2}-96p_{3})n.
\end{eqnarray*}
Now, we proceed by showing the following results.

\begin{thm}
Suppose Hypotheses 1 and 2 are true, then the next results hold. $(i)$ The variance of the Schultz index is denoted by
\begin{eqnarray*}
Var(S(COC_{n}))&=&\frac{1}{30}\bigg({\sigma}_{2}^{2}n^{5}-5r_{2}n^{4}+10\tilde{\sigma}^{2}_{2}n^{3}+(65r_{2}-30{\sigma}^{2}_{2}-45\tilde{\sigma}^{2}_{2})n^{2}\\
&&+(59{\sigma}^{2}_{2}+65\tilde{\sigma}^{2}_{2}-120r_{2})n+(60r_{2}-30{\sigma}^{2}_{2}-30\tilde{\sigma}^{2}_{2})\bigg),
\end{eqnarray*}
where
\begin{align*}
\sigma^{2}_{2}&=576^{2}p_{1}+864^{2}p_{2}+1152^{2}p_{3}+1440^{2}(1-p_{1}-p_{2}-p_{3})\\
&~~~~-\big(576p_{1}+864p_{2}+1152p_{3}+1440(1-p_{1}-p_{2}-p_{3})\big)^{2},\\
\tilde{\sigma}_{2}^{2}&=40^{2}p_{1}+328^{2}p_{2}+616^{2}p_{3}+904^{2}(1-p_{1}-p_{2}-p_{3})\\
&~~~~-\big(40p_{1}+328p_{2}+616p_{3}+904(1-p_{1}-p_{2}-p_{3})\big)^{2},\\
r_{2}&=-\big(576p_{1}+864p_{2}+1152p_{3}+1440(1-p_{1}-p_{2}-p_{3})\big)\big(40p_{1}+328p_{2}+616p_{3}+904(1-p_{1}-p_{2}-p_{3})\big)\\
&~~~~+576\cdot 40p_{1}+864\cdot 328p_{2}+1152\cdot 616p_{3}+1440\cdot 904(1-p_{1}-p_{2}-p_{3}).
\end{align*}
(ii)  For $n\rightarrow \infty$, $S(COC_{n})$ asymptotically obeys normal distributions. One has
\begin{eqnarray*}
\lim_{n\rightarrow \propto}\sup_{a\in \mathbb{R}}\mid \mathbb{P}\big(\frac{S(COC_{n})-\mathbb{E}(S(COC_{n}))}{\sqrt{Var(S(COC_{n}))}}\leq a\big)-{\int}^{a}_{-\infty}\frac{1}{\sqrt{2\pi}}e^{-\frac{t^{2}}{2}}dt \mid =0.
\end{eqnarray*}
\end{thm}

\noindent\textbf{Proof of Theorem 3.2.} Let $D_{n}=\sum\limits_{v \in V_{COC_{n}}}(8deg(v)+18)dis(u_{n},v)$. Then by Eq. (2.13), we obtain
\begin{eqnarray}
S(COC_{n+1})=S(COC_{n})+D_{n}+824n+248.
\end{eqnarray}
According to the previous proof of Theorem 3.1, the four Facts are obtained.

\textbf{Fact 3.2.1.} $D_{n}Z_{n}^{1}=(D_{n-1}+576n-40)Z_{n}^{1}.$

\noindent\textbf{Proof.} If $Z_{n}^{1}$ = 0, the above result is distinct. So we take into account $Z_{n}^{1}$ = 1, which indicates $COC_{n} \rightarrow COC^{1}_{n+1}$. In this case, $u_{n}$ (of $COC_{n}$) coincides with $x_{2}$ or $x_{8}$ (of $O_{n}$), see Fig. 3. In this situation, by using Eq. (2.7)-Eq. (2.8), $D_{n}$ becomes
\begin{align*}
&\sum\limits_{v \in V_{COC_{n}}}(8deg(v)+18)dis(x_{2},v)\\
&=\sum_{v \in V_{COC_{n-1}}}(8deg(v)+18)dis(x_{2},v)+\sum_{v \in V_{O_{n}}}(8deg(v)+18)dis(x_{2},v)\\
&=\sum_{v \in V_{COC_{n-1}}}(8deg(v)+18)(dis(u_{n-1},v)+2)+8\times 33+18\times 16\\
&=\sum_{v \in V_{COC_{n-1}}}(8deg(v)+18)dis(u_{n-1},v)+2\cdot \big(18\times8(n-1)+8\times(18n-19)\big)+552\\
&=D_{n-1}+576n-40.
\end{align*}
Thus, we conclude the desired Fact.

\textbf{Fact 3.2.2.} $D_{n}Z_{n}^{2}=(D_{n-1}+864n-328)Z_{n}^{2}.$

As that in the proof of Fact 3.2.1, we only consider the fact $Z_{n}^{2}=1$, that is $COC_{n} \rightarrow COC^{2}_{n+1}$. The proof is similar and details are omitted.

\textbf{Fact 3.2.3.} $D_{n}Z_{n}^{3}=(D_{n-1}+1152n-616)Z_{n}^{3}.$

Similarly, we only consider the fact $Z_{n}^{3}=1$, that is $COC_{n} \rightarrow COC^{3}_{n+1}$. We omit the details.

\textbf{Fact 3.2.4.} $ D_{n}Z_{n}^{4}=(D_{n-1}+1440n-904)Z_{n}^{4}.$

At the same way, we just consider the fact $Z_{n}^{4}=1$, that is $COC_{n} \rightarrow COC^{4}_{n+1}$. The details are omitted here.

Noting that $Z_{n}^{1}+Z_{n}^{2}+Z_{n}^{3}+Z_{n}^{4}=1$, by the above discussions, it holds that
\begin{eqnarray*}
D_{n}&=&D_{n}(Z_{n}^{1}+Z_{n}^{2}+Z_{n}^{3}+Z_{n}^{4})\\
&=&(D_{n-1}+576n-40)Z_{n}^{1}+(D_{n-1}+864n-328)Z_{n}^{2}+(D_{n-1}+1152n-616)Z_{n}^{3}\\
&&+(D_{n-1}+1440n-904)Z_{n}^{4}\\
&=&D_{n-1}+(576Z_{n}^{1}+864Z_{n}^{2}+1152Z_{n}^{3}+1440Z_{n}^{4})n-(40Z_{n}^{1}+328Z_{n}^{2}+616Z_{n}^{3}+904Z_{n}^{4})\\
&=&D_{n-1}+nU^{'}_{n}-V^{'}_{n},
\end{eqnarray*}
where for each $n$,

~~~~~~~~$U^{'}_{n}=576Z_{n}^{1}+864Z_{n}^{2}+1152Z_{n}^{3}+1440Z_{n}^{4}$,~~~~$V^{'}_{n}=40Z_{n}^{1}+328Z_{n}^{2}+616Z_{n}^{3}+904Z_{n}^{4}$.
\\\text{Therefore, by Eq. (3.19),}\par

\begin{align}
S(COC_{n})=& S(COC_{1})+\sum_{a=1}^{n-1}D_{a}+\sum_{a=1}^{n-1}(824a+248)\nonumber\\
=& S(COC_{1})+\sum_{a=1}^{n-1}\big(\sum_{b=1}^{a-1}(D_{b+1}-D_{b})+D_{1}\big)+\sum_{a=1}^{n-1}(824a+248)\nonumber\\
=& S(COC_{1})+\sum_{a=1}^{n-1}\sum_{b=1}^{a-1}(D_{b+1}-D_{b})+(n-1)D_{1}+\sum_{a=1}^{n-1}(824a+248)\nonumber\\
=& S(COC_{1})+\sum_{a=1}^{n-1}\sum_{b=1}^{a-1}\big((b+1)U^{'}_{b+1}-V^{'}_{b+1}\big)+O(n^2)\nonumber.
\end{align}

Suppose that

~~~~~~~~~~~~~~~~~~~~~~$Var(U^{'}_{b})=\sigma^{2}_{2}$,~~~~$Var(V^{'}_{b})=\tilde{\sigma}^{2}_{2}$,~~~~$Cov(U^{'}_{b},V^{'}_{b})=r_{2}$.

If we substitute $Gut(COC_{n})$ by $S(COC_{n})$ in the proof of Theorem 3.1, the rest of the proof is identical to the proof of Theorem 3.1 and the details are omitted here.\hfill\rule{1ex}{1ex}\

\section{The limiting behaviours for $Kf^{*}(COC_{n})$ and $Kf^{+}(COC_{n})$ of a random cyclooctane chain}
 \ \ \ \ \
Following the previous sections, we determine the limiting behaviours for $Kf^{*}(COC_{n})$ and $Kf^{+}(COC_{n})$. For $n\geq 2$, $COC_{n+1}$ is organized by adding a new terminal octagon $O_{n+1}$ to $COC_{n}$ by an edge, where $O_{n+1}$ are labelled as $x_{1},~ x_{2}, ~x_{3}, ~x_{4},~ x_{5},~ x_{6},~ x_{7},~ x_{8}$. For all $v \in V_{COC_{n}} \subseteq V_{COC_{n+1}}$, one has
\begin{equation}
\begin{split}
&r(x_{1},v)=r(u_{n},v)+1,~~~~~~~~~~r(x_{2},v)=r(u_{n},v)+1+\frac{7}{8},~~~r(x_{3},v)=r(u_{n},v)+1+\frac{3}{2},\\
&r(x_{4},v)=r(u_{n},v)+1+\frac{15}{8},~~~r(x_{5},v)=r(u_{n},v)+1+2,~~~r(x_{6},v)=r(u_{n},v)+1+\frac{15}{8},\\
&r(x_{7},v)=r(u_{n},v)+1+\frac{3}{2},~~~~~r(x_{8},v)=r(u_{n},v)+1+\frac{7}{8},\\
&\sum_{v \in V_{COC_{n-1}}}deg_{COC_{n}}(v)=18n-19, ~~and ~~\sum_{v \in V_{COC_{n}}}deg_{COC_{n+1}}(v)=18n-1.
\end{split}
\end{equation}
Meanwhile,
\begin{equation}
\begin{split}
&\sum_{i=1}^{8}deg(x_{i})r(x_{1},x_{i})=21,~~~\sum_{i=1}^{8}deg(x_{i})r(x_{2},x_{i})=\frac{175}{8},~~~\sum_{i=1}^{8}deg(x_{i})r(x_{3},x_{i})=\frac{45}{2},\\
&\sum_{i=1}^{8}deg(x_{i})r(x_{4},x_{i})=\frac{183}{8},~~~\sum_{i=1}^{8}deg(x_{i})r(x_{5},x_{i})=23,~~~\sum_{i=1}^{8}deg(x_{i})r(x_{6},x_{i})=\frac{183}{8},\\
&\sum_{i=1}^{8}deg(x_{i})r(x_{7},x_{i})=\frac{45}{2},~~~\sum_{i=1}^{8}deg(x_{i})r(x_{8},x_{i})=\frac{175}{8}.\\
\end{split}
\end{equation}
Together with Eq. (1.4), in [\cite{a29}, Theorem 3], the paper showed that
\begin{eqnarray*}
\mathbb{E}(Kf^{*}(COC_{n}))&=&6(27-\frac{81}{8}p_{1}-\frac{9}{2}p_{2}-\frac{9}{8}p_{3})n^{3}+18(2+\frac{81}{8}p_{1}+\frac{9}{2}p_{2}+\frac{9}{8}p_{3})n^{2}\\
&&-12(\frac{29}{12}+\frac{81}{8}p_{1}+\frac{9}{2}p_{2}+\frac{9}{8}p_{3})n-1.
\end{eqnarray*}
Then, we go on to show the following results.

\begin{thm}
Suppose Hypotheses 1 and 2 are true, there are the next main results. $(i)$ The variance of $Kf^{*}(COC_{n})$ is denoted by
\begin{eqnarray*}
Var(Kf^{*}(COC_{n}))&=&\frac{1}{30}\bigg({\sigma}^{2}_{3}n^{5}-5r_{3}n^{4}+10\tilde{\sigma}^{2}_{3}n^{3}+(65r_{3}-30{\sigma}^{2}_{3}-45\tilde{\sigma}^{2}_{3})n^{2}\\
&&+(59{\sigma}^{2}_{3}+65\tilde{\sigma}^{2}_{3}-120r_{3})n+(60r_{3}-30{\sigma}^{2}_{3}-30\tilde{\sigma}^{2}_{3})\bigg),
\end{eqnarray*}
where
\begin{align*}
\sigma^{2}_{3}&=(\frac{1215}{2})^{2}p_{1}+810^{2}p_{2}+(\frac{1863}{2})^{2}p_{3}+972^{2}(1-p_{1}-p_{2}-p_{3})\\
&~~~~-\big(\frac{1215}{2}p_{1}+810p_{2}+\frac{1863}{2}p_{3}+972(1-p_{1}-p_{2}-p_{3})\big)^{2},\\
\tilde{\sigma}^{2}_{3}&=(\frac{495}{2})^{2}p_{1}+450^{2}p_{2}+(\frac{1143}{2})^{2}p_{3}+612^{2}(1-p_{1}-p_{2}-p_{3})\\
&~~~~-\big(\frac{495}{2}p_{1}+450p_{2}+\frac{1143}{2}p_{3}+612(1-p_{1}-p_{2}-p_{3})\big)^{2},\\
r_{3}&=-\big(\frac{1215}{2}p_{1}+810p_{2}+\frac{1863}{2}p_{3}+972(1-p_{1}-p_{2}-p_{3})\big)\big(\frac{495}{2}p_{1}+450p_{2}+\frac{1143}{2}p_{3}+612(1-p_{1}-p_{2}-p_{3})\big)\\
&~~~~+\frac{1215}{2}\cdot \frac{495}{2}p_{1}+810\cdot 450p_{2}+\frac{1863}{2}\cdot \frac{1143}{2}p_{3}+972\cdot 612(1-p_{1}-p_{2}-p_{3}).
\end{align*}
(ii)  For $n\rightarrow \infty$, $Kf^{*}(COC_{n})$ asymptotically obeys normal distributions. One has
\begin{eqnarray*}
\lim_{n\rightarrow \propto}\sup_{a\in \mathbb{R}}\mid \mathbb{P}\big(\frac{Kf^{*}(COC_{n})-\mathbb{E}(Kf^{*}(COC_{n}))}{\sqrt{Var(Kf^{*}(COC_{n}))}}\leq a\big)-{\int}^{a}_{-\infty}\frac{1}{\sqrt{2\pi}}e^{-\frac{t^{2}}{2}}dt \mid =0.
\end{eqnarray*}
\end{thm}

\noindent\textbf{\bf Proof of Theorem 4.1.} Subsitituting $n+1$ into [\cite{a29}, Eq. (11)],  we see
\begin{eqnarray*}
Kf^{*}(COC_{n+1})=Kf^{*}(COC_{n})+18\sum_{v\in V_{COC_{n}}}deg(v)r(u_{n},v)+684n+151.
\end{eqnarray*}
Let $P_{n}=18\sum\limits_{v \in V_{COC_{n}}}deg(v)r(u_{n},v)$, we have
\begin{eqnarray}
Kf^{*}(COC_{n+1})=Kf^{*}(COC_{n})+P_{n}+684n+151.
\end{eqnarray}

Recalling that $Z_{n}^{1}$, $Z_{n}^{2}$, $Z_{n}^{3}$ and $Z_{n}^{4}$ are random variables which show the way to construct $COC_{n+1}$ from $COC_{n}$. We get the following four Facts.

\noindent\textbf{\bf Fact 4.1.1.} $P_{n}Z_{n}^{1}=(P_{n-1}+\frac{1215}{2}n-\frac{495}{2})Z_{n}^{1}.$

\noindent\textbf{\bf Proof.} If $Z_{n}^{1}$ = 0, the result is obvious. So we just take into account $Z_{n}^{1}$ = 1, which implies $COC_{n} \rightarrow COC^{1}_{n+1}$. In this case, $u_{n}$ (of $COC_{n}$) overlaps with $x_{2}$ or $x_{8}$ (of $O_{n}$), see Fig. 3. In this situation, by using Eq. (4.20)-Eq. (4.21), $P_{n}$ becomes
\begin{eqnarray*}
18\sum_{v \in V_{COC_{n}}}d(v)r(x_{2},v)&=&18\sum_{v \in V_{COC_{n-1}}}deg(v)r(x_{2},v)+18\sum_{v \in V_{O_{n}}}deg(v)r(x_{2},v)\\
&=&18\sum_{v \in V_{COC_{n-1}}}deg(v)(r(u_{n-1},v)+1+\frac{7}{8})+18\times \frac{175}{8}\\
&=&P_{n-1}+18\times \frac{15}{8}\times (18n-19)+18\times \frac{175}{8}\\
&=&P_{n-1}+\frac{1215}{2}n-\frac{495}{2}.
\end{eqnarray*}
Thus, we conclude the desired Fact.

\noindent\textbf{\bf Fact 4.1.2.} $P_{n}Z_{n}^{2}=(P_{n-1}+810n-450)Z_{n}^{2}.$

As that in the proof of Fact 4.1.1, we only consider the fact $Z_{n}^{2}=1$, that is $COC_{n} \rightarrow COC^{2}_{n+1}$. The proof is similar and details are omitted.

\noindent\textbf{\bf Fact 4.1.3.} $P_{n}Z_{n}^{3}=(P_{n-1}+\frac{1863}{2}n-\frac{1143}{2})Z_{n}^{3}.$

Similarly, we only consider the fact $Z_{n}^{3}=1$, that is $COC_{n} \rightarrow COC^{3}_{n+1}$. At the same way, we omit the details.

\noindent\textbf{\bf Fact 4.1.4.} $P_{n}Z_{n}^{4}=(P_{n-1}+972n-612)Z_{n}^{4}.$

We only consider the fact $Z_{n}^{4}=1$, that is $COC_{n} \rightarrow COC^{4}_{n+1}$. details are omitted here.

Noting that $Z_{n}^{1}+Z_{n}^{2}+Z_{n}^{3}+Z_{n}^{4}=1$, by the above discussions, it holds that
\begin{eqnarray*}
P_{n}&=&P_{n}(Z_{n}^{1}+Z_{n}^{2}+Z_{n}^{3}+Z_{n}^{4})\\
&=&(P_{n-1}+\frac{1215}{2}n-\frac{495}{2})Z_{n}^{1}+(P_{n-1}+810n-450)Z_{n}^{2}+(P_{n-1}+\frac{1863}{2}n-\frac{1143}{2})Z_{n}^{3}\\
&&+(P_{n-1}+972n-612)Z_{n}^{4}\\
&=&P_{n-1}+(\frac{1215}{2}Z_{n}^{1}+810Z_{n}^{2}+\frac{1863}{2}Z_{n}^{3}+972Z_{n}^{4})n-(\frac{495}{2}Z_{n}^{1}+450Z_{n}^{2}+\frac{1143}{2}Z_{n}^{3}+612Z_{n}^{4})\\
&=&P_{n-1}+n\hat{U}_{n}-\hat{V}_{n},
\end{eqnarray*}
where for each $n$,

~~~~~~~$\hat{U}_{n}=\frac{1215}{2}Z_{n}^{1}+810Z_{n}^{2}+\frac{1863}{2}Z_{n}^{3}+972Z_{n}^{4}$,~~~~$\hat{V}_{n}=\frac{495}{2}Z_{n}^{1}+450Z_{n}^{2}+\frac{1143}{2}Z_{n}^{3}+612Z_{n}^{4}$.
\\\text{Therefore, by Eq. (4.22),}\par
\begin{align}
Kf^{*}(COC_{n})=& Kf^{*}(COC_{1})+\sum_{a=1}^{n-1}P_{a}+\sum_{a=1}^{n-1}(684a+151)\nonumber\\
=& Kf^{*}(COC_{1})+\sum_{a=1}^{n-1}(\sum_{b=1}^{a-1}(P_{b+1}-P_{b})+P_{1})+\sum_{a=1}^{n-1}(684a+151)\nonumber\\
=& Kf^{*}(COC_{1})+\sum_{a=1}^{n-1}\sum_{b=1}^{a-1}((b+1)\hat{U}_{b+1}-\hat{V}_{b+1})+(n-1)P_{1}+\sum_{a=1}^{n-1}(684a+151)\nonumber\\
=& Kf^{*}(COC_{1})+\sum_{a=1}^{n-1}\sum_{b=1}^{a-1}((b+1)\hat{U}_{b+1}-\hat{V}_{b+1})+O(n^2).
\end{align}

By direct calculation, we put

~~~~~~~~~~~~~~~~~~~~$Var(\hat{U}_{b})=\sigma^{2}_{3}$,~~~~$Var(\hat{V}_{b})=\tilde{\sigma}^{2}_{3}$,~~~~$Cov(\hat{U}_{b},\hat{V}_{b})=r_{3}$.

By the properties of variance, Eq. (4.23) and interchanging the orders of $a$ and $b$, it follows that
\begin{align*}
&Var(Kf^{*}(COC_{n}))\\
&=Var\bigg(\sum_{a=1}^{n-1}\sum_{b=1}^{a-1}\big((b+1)\hat{U}_{b+1}-\hat{V}_{b+1}\big)\bigg)=Var\bigg(\sum_{b=1}^{n-2}\sum_{a=b+1}^{n-1}\big((b+1)\hat{U}_{b+1}-\hat{V}_{b+1}\big)\bigg)\\
&=\sum_{b=1}^{n-2}(n-b-1)^{2}Cov\big((b+1)\hat{U}_{b+1}-\hat{V}_{b+1},(b+1)\hat{U}_{b+1}-\hat{V}_{b+1}\big)\\
&=\sum_{b=1}^{n-2}(n-b-1)^{2}\big((b+1)^2\sigma^{2}_{3}-2(b+1)r_{3}+\tilde{\sigma}^{2}_{3}\big).\\
\end{align*}
The above formula indicates the result $Theorem ~4.1.~ (i)$ with the help of a computer.

The rest is the same as that in the proof of $Theorem ~3.1.~~ (ii)$ and details are omitted here.\hfill\rule{1ex}{1ex}\

Similarly, together with Eq. (1.5), in [\cite{a29}, Theorem 5], the paper put forward that
\begin{eqnarray*}
\mathbb{E}(Kf^{+}(COC_{n}))&=&6(24-9p_{1}-4p_{2}-p_{3})n^{3}+18(\frac{61}{18}+9p_{1}+4p_{2}+p_{3})n^{2}\\
&&-12(\frac{37}{12}+9p_{1}+4p_{2}+p_{3})n.
\end{eqnarray*}
We proceed by showing the following result about the expatiatory formula of the variance of $Kf^{+}(COC_{n})$.

\begin{thm}
 Suppose Hypotheses 1 and 2 are true, there are the next main results. $(i)$ The variance of $Kf^{+}(COC_{n})$ is denoted by
\begin{eqnarray*}
Var(Kf^{+}(COC_{n}))&=&\frac{1}{30}\bigg({\sigma}^{2}_{4}n^{5}-5r_{4}n^{4}+10\tilde{\sigma}^{2}_{4}n^{3}+(65r_{4}-30{\sigma}^{2}_{4}-45\tilde{\sigma}^{2}_{4})n^{2}\\
&&+(59{\sigma}^{2}_{4}+65\tilde{\sigma}^{2}_{4}-120r_{4})n+(60r_{4}-30{\sigma}^{2}_{4}-30\tilde{\sigma}^{2}_{4})\bigg),
\end{eqnarray*}
where
\begin{align*}
\sigma^{2}_{4}&=540^{2}p_{1}+720^{2}p_{2}+828^{2}p_{3}+864^{2}(1-p_{1}-p_{2}-p_{3})\\
&~~~~-\big(540p_{1}+720p_{2}+828p_{3}+864(1-p_{1}-p_{2}-p_{3})\big)^{2},\\
\tilde{\sigma}^{2}_{4}&=191^{2}p_{1}+371^{2}p_{2}+479^{2}p_{3}+515^{2}(1-p_{1}-p_{2}-p_{3})\\
&~~~~-\big(191p_{1}+371p_{2}+479p_{3}+515(1-p_{1}-p_{2}-p_{3})\big)^{2},\\
r_{4}&=-\big(540p_{1}+720p_{2}+828p_{3}+864(1-p_{1}-p_{2}-p_{3})\big)\big(191p_{1}+371p_{2}+479p_{3}+515(1-p_{1}-p_{2}-p_{3})\big)\\
&~~~~+540\cdot 191p_{1}+720\cdot 371p_{2}+828\cdot 479p_{3}+864\cdot 515(1-p_{1}-p_{2}-p_{3}).
\end{align*}
(ii) For $n \rightarrow \infty$, $Kf^{+}(COC_{n})$ asymptotically obeys normal distributions. One has
\begin{eqnarray*}
\lim_{n\rightarrow \propto}\sup_{a\in \mathbb{R}}\mid \mathbb{P}\big(\frac{Kf^{+}(COC_{n})-\mathbb{E}(Kf^{+}(COC_{n}))}{\sqrt{Var(Kf^{+}(COC_{n}))}}\leq a\big)-{\int}^{a}_{-\infty}\frac{1}{\sqrt{2\pi}}e^{-\frac{t^{2}}{2}}dt \mid =0.
\end{eqnarray*}
\end{thm}

\noindent\textbf{\bf Proof of Theorem 4.2.}  Subsitituting $n+1$ into [\cite{a29}, Eq. (15)] yields, one sees that
\begin{eqnarray*}
Kf^{+}(COC_{n+1})=Kf^{+}(COC_{n})+\sum_{v\in V_{COC_{n}}}(18+8deg(v))r(u_{n},v)+637n+160.
\end{eqnarray*}
 Let $Q_{n}=\sum\limits_{v \in V_{COC_{n}}}(18+8deg(v))r(u_{n},v)$. Then, we obtain
\begin{eqnarray}
Kf^{+}(COC_{n+1})=Kf^{+}(COC_{n})+Q_{n}+637n+160.
\end{eqnarray}

Recalling that $Z_{n}^{1}$, $Z_{n}^{2}$, $Z_{n}^{3}$ and $Z_{n}^{4}$ are random variables, we obtain the following Facts.

\noindent\textbf{\bf Fact 4.2.1.} $Q_{n}Z_{n}^{1}=(Q_{n-1}+540n-191)Z_{n}^{1}.$

\noindent\textbf{\bf Proof.} If $Z_{n}^{1}$ = 0, the result is distinct. Then, we just take into account $Z_{n}^{1}$ = 1, which indicates $COC_{n} \rightarrow COC^{1}_{n+1}$. In this fact, $u_{n}$ (of $COC_{n}$) coincides with $x_{2}$ or $x_{8}$ (of $O_{n}$), see Fig. 3. In this situation, $Q_{n}$ becomes
\begin{align*}
\begin{split}
&\sum_{v \in V_{COC_{n}}}(18+8deg(v))r(x_{2},v)\\
&=\sum_{v \in V_{COC_{n-1}}}(18+8deg(v))r(x_{2},v)+\sum_{v \in V_{O_{n}}}(18+8deg(v))r(x_{2},v)\\
&=\sum_{v \in V_{COC_{n-1}}}(18+8deg(v))(r(u_{n-1},v)+1+\frac{7}{8})+18\sum_{v \in V_{O_{n}}}r(x_{2},v)+8\sum_{v \in V_{O_{n}}}deg(v)r(x_{2},v)\\
&=\sum_{v \in V_{COC_{n-1}}}(18+8deg(v))r(u_{n-1},v)+\frac{15}{8}\sum_{v \in V_{COC_{n-1}}}(18+8deg(v))+8\times \frac{175}{8}+18\times \frac{84}{8}\\
&=Q_{n-1}+540n-191.
\end{split}
\end{align*}
Thus, we obtain the desired Fact.

\noindent\textbf{\bf Fact 4.2.2.} $Q_{n}Z_{n}^{2}=(Q_{n-1}+720n-371)Z_{n}^{2}.$

Similar to the proof of Fact 4.2.1, we only consider the fact $Z_{n}^{2}=1$, that is $COC_{n} \rightarrow COC^{2}_{n+1}$. In the same way, we omit the details.

\noindent\textbf{\bf Fact 4.2.3.} $Q_{n}Z_{n}^{3}=(Q_{n-1}+828n-479)Z_{n}^{3}.$

Similarly, we only consider the fact $Z_{n}^{3}=1$, that is $COC_{n} \rightarrow COC^{3}_{n+1}$. The proof is also similar to the above facts and we omit the details.

\noindent\textbf{\bf Fact 4.2.4.} $Q_{n}Z_{n}^{4}=(Q_{n-1}+864n-515)Z_{n}^{4}.$

Considering the fact $Z_{n}^{4}=1$ which is $COC_{n} \rightarrow COC^{4}_{n+1}$, the details are omitted here..

Noting that $Z_{n}^{1}+Z_{n}^{2}+Z_{n}^{3}+Z_{n}^{4}=1$, by the above discussions, it holds that
\begin{eqnarray*}
Q_{n}&=&Q_{n-1}+(540Z_{n}^{1}+720Z_{n}^{2}+828Z_{n}^{3}+864Z_{n}^{4})n-(191Z_{n}^{1}+371Z_{n}^{2}+479Z_{n}^{3}+515Z_{n}^{4})\\
&=&Q_{n-1}+n\tilde{U}_{n}-\tilde{V}_{n},
\end{eqnarray*}
where for each $n$,

~~~~$\tilde{U}_{n}=540Z_{n}^{1}+720Z_{n}^{2}+828Z_{n}^{3}+864Z_{n}^{4}$,~~~~$\tilde{V}_{n}=191Z_{n}^{1}+371Z_{n}^{2}+479Z_{n}^{3}+515Z_{n}^{4}$.
\\\text{Therefore, by Eq. (4.24),}\par

\begin{align}
Kf^{+}(COC_{n})=& Kf^{+}(COC_{1})+\sum_{a=1}^{n-1}Q_{a}+\sum_{a=1}^{n-1}(637a+160)\nonumber\\
=& Kf^{+}(COC_{1})+\sum_{a=1}^{n-1}\big(\sum_{b=1}^{a-1}(Q_{b+1}-Q_{b})+Q_{1}\big)+\sum_{a=1}^{n-1}(637a+160)\nonumber\\
=& Kf^{+}(COC_{1})+\sum_{a=1}^{n-1}\sum_{b=1}^{a-1}\big((b+1)\tilde{U}_{b+1}-\tilde{V}_{b+1}\big)+(n-1)Q_{1}+\sum_{a=1}^{n-1}(637a+160)\nonumber\\
=& Kf^{+}(COC_{1})+\sum_{a=1}^{n-1}\sum_{b=1}^{a-1}\big((b+1)\tilde{U}_{b+1}-\tilde{V}_{b+1}\big)+O(n^2)\nonumber.
\end{align}

Suppose that

~~~~~~~~~~$Var(\tilde{U}_{b})=\sigma^{2}_{4}$,~~~~$Var(\tilde{V}_{b})=\tilde{\sigma}^{2}_{4}$,~~~~$Cov(\tilde{U}_{b},\tilde{V}_{b})=r_{4}$.

The rest of proof is similar to the above Theorem, and the details are omitted here.\hfill\rule{1ex}{1ex}\

\section*{Acknowledgments}

The authors would like to express their sincere gratitude to the editor and anonymous referees for valuable suggestions, which led to great deal of improvement of the original paper.

\section*{Disclosure statement}

No potential conflict of interest was reported by the authors.

\section*{Funding}

This work was supported in part by Anhui Provincial Natural Science Foundation under Grant 2008085J01 and Natural Science Fund of Education Department of Anhui Province under Grant KJ2020-A0478.

\end{document}